\documentclass{article}
\usepackage{graphicx}
\usepackage{amsmath, amsthm, amssymb, amscd}
\usepackage{cite}

\def\beq{\begin{equation}}
\def\eeq{\end{equation}}

\begin{document}

\title{Geodesics on Surfaces with Helical Symmetry: 
Cavatappi Geometry
}
\author{
Robert T. Jantzen\\ 
Department of Mathematics and Statistics\\
Villanova University\\
Villanova, PA 19085 USA
}

\date{December 31, 2012}
\maketitle

\begin{abstract}
A 3-parameter family of helical tubular surfaces obtained by screw revolving a circle provides a useful pedagogical example of how to study geodesics on a surface that admits a 1-parameter symmetry group, but is not as simple as a surface of revolution like the torus which it contains as a special case. It serves as a simple example of helically symmetric surfaces which are the generalizations of surfaces of revolution in which an initial plane curve is screw-revolved around an axis in its plane. The physics description of geodesic motion on these surfaces requires a slightly more involved effective potential approach than the torus case due to the nonorthogonal coordinate grid necessary to describe this problem. Amazingly this discussion allows one to very nicely describe the geodesics of the surface of the more complicated ridged cavatappi pasta.

\end{abstract}

\newpage
\begin{figure}[t] 
\typeout{*** EPS figure cavatappo}
\vglue-2.6cm
\begin{center}
\includegraphics[scale=0.7]{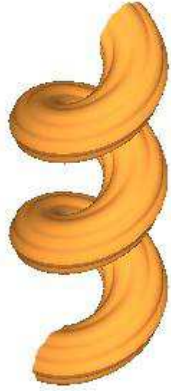}
\end{center}\vglue-2cm
\caption{The ridged cavatappo (singular) / cavatappi (plural) pasta shape.}

\label{fig:cavatappo}
\end{figure}

\section{Introduction}

In looking for examples of interesting surfaces to use in playing with geodesics, one starts with the plane and cylinder (flat, rotationally and translation symmetric), moves on to the sphere (rotationally symmetric), whose great circle geometry was long used in celestial navigation even before non-Euclidean geometry and then differential geometry found them to be nice mathematical objects of study. These are all maximally symmetric (both isotropic at every point: every direction is equivalent, and homogeneous: each point is equivalent to every other point) with 3-dimensional symmetry groups. After this, surfaces of revolution (which have rotational symmetry about an axis) present themselves as candidates for study, with 
the 2-parameter family of tori as the first inhomogeneous but still rotationally symmetric simple surface one encounters next yet still constructed from circles, demonstrating both positive and negative curvature, and containing the sphere as a special case. This provides a nice test case to study geodesics using the physics approach to motion along curves constrained to lie on a surface, revealing lots of fun mathematics developed in detail in a previous investigation  \cite{bobtorus}.

What next? It turns out that the tubular surface based on a helix is a 3-parameter family which contains the torus as a special case where the helix degenerates to a circle, but it still has a 1-parameter symmetry group of corkscrew rotations (a linked rotation about and translation along the central axis of the helix) that makes its study reducible to motion in 1-dimension, like the torus. In fact amazingly many of the properties of geodesic motion on the torus apply to the helical tube. One can imagine constructing one screw-revolution of this tube  by cutting a horizontal torus by a half-plane through its vertical axis of symmetry along a vertical circular cross-section (a meridian), and then stretching the two ends vertically in opposite directions until the circle of centers of all such cross-sections stretches into a segment of a helix about the axis of symmetry. One can then extend the terminology of inner and outer (helical) equators, (circular) meridians, (helical) parallels, and northern and southern polar helices in a natural way from the corresponding curves on the torus, in turn taken from the sphere. 
Remarkably the equators remain geodesics of the surface as in the torus case, while 
the study of general geodesics can be reduced to the 1-dimensional problem of ``radial" motion along the meridians due to the conserved ``screw-angular momentum" associated with the 1-parameter symmetry group, generalizing the conserved vertical component of angular momentum of the corresponding torus case. The terminology of bound and unbound geodesics also generalizes to the new case, with a qualitative correspondence between the bound geodesics oscillating about the outer equator but not crossing the inner equator, and those unbound geodesics which do cross the inner equator.

This larger class of surfaces containing the tori is useful in that it  adds a new element to the discussion: a nonorthogonal natural coordinate grid on the surface. This is easily orthogonalized in 2-dimensions by completing the square to adapt a new basis of the tangent spaces to the orthogonal decomposition with respect to the parallels  which are the helical  orbits of the 1-dimensional symmetry group of this geometry. Once done, all of the qualitative features of geodesic motion on tori generalize to the larger family of helical tubular orbits. We are all familiar with this symmetry: turn a metal screw in its hole by 360 degrees and it returns to its original orientation except moved farther into the hole. The angle of rotation is proportional to a translation along the axis of the rotation.
 The present discussion is understood to be a sequel to the torus article \cite{bobtorus}, showing how to take into account the nonorthogonal coordinate system to reduce the geodesic discussion to one qualitatively like the simpler torus case.

For particular values of the shape parameters of our family of helical tubular surfaces, two and a half turns of the tube models extremely well the ridgeless (liscio) cavatappi pasta shape (also known as gobetti, spirali or cellentani), which is a common corkscrew version of the macaroni family of tube pastas whose design aspects were recently documented by George Legendre in this elegant book \textit{Pasta by Design} \cite{pastadesign,blog}. The default ``rigato" (ridged) cavatappo that is actually found in markets is a slight generalization of the helical tubular shape within the class of surfaces admitting a helical symmetry and although its geodesic equations require several pages to write out, they are easily handled numerically by the computer algebra system required to derive them. For this pasta shape, this exercise takes the mathematics of pasta one step beyond the design aspects documented by the architect George Legendre to the dynamics of its intrinsic geometry as approached by a physicist.

Although the screw-rotation symmetry is familiar to most of us through its everyday representation in bolts and metal sheet screws (but not wood screws which taper off) and cork screw wine bottle openers, the corresponding mathematical symmetry group is surprisingly absent as an example in mathematical discussions of continuous transformation groups \cite{wikiscrew} or in toy applications in classical mechanical systems. However, the Archimedes screw has been with us for millenia and still finds useful applications in large scale industrial screws that benefit society in various ways \cite{archimedes}. Surprisingly, a careful mathematical study of improving the traditional design was only done recently by Chris Rorres with the help of a computer algebra system \cite{rorres}, with a followup study resulting in an even more optimal design.

\begin{figure}
\typeout{*** EPS figure torus1}
\vglue0cm
\begin{center}
 \includegraphics[scale=0.3]{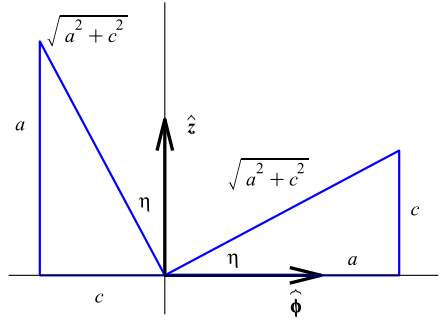}
\end{center}
\vglue0cm
\caption{
An illustration of the inclination angle $\eta$  from the horizontal of the central helix of the surface in $\phi$-$z$ plane of the tangent space at a point along that curve, whose unit tangent vector is aligned with the hypotenuse of the right right triangle, and binormal with the hypotenuse of the left right triangle, while the normal $-\hat \rho$ points along the cylindrical radial coordinate line towards the symmetry axis, where $(\hat \rho, \hat\phi,\hat z)$ are unit vectors along the cylindrical coordinate lines.
The value of the inclination angle  $\eta=\arctan(8/15)\approx 28.1^\circ$ illustrated here is very nearly equal to Legendre's choice for his cavatappo design parameter but corresponds to the fourth Pythagorean triple $(8,15,17)$ right triangle leading to rational values of the surface parameters: $(a,b,c)=\left( \frac32,1,\frac45 \right)=\frac32\left( 1,\frac23,\frac8{15} \right)$.
} 
\label{fig:inclination}
\end{figure}

\section{Helical tubular surfaces}

Starting from the usual Cartesian coordinates $(x,y,z)$ in Euclidean space where the element of differential arclength $ds$ satisfies an iteration of the Pythagorean theorem  $ds^2=dx^2+dy^2+dz^2$ (the ``line element" representation of the ``metric" of Euclidean space), a simple tubular surface built from a helix can be represented by the following parametrization \cite{pastadesign,blog}
\beq\label{eq:helix}
  x = (a+b \cos \chi) \cos \phi\,,\
  y = (a+b \cos \chi) \sin \phi\,,\
  z = c\, \phi+b \sin \chi \,,
\eeq
or equivalently as $\rho=a+b\, \cos \chi, z = c\, \phi+b\, \sin \chi$ in the usual cylindrical coordinates $(\rho,\phi,z)$ related to the Cartesian ones by 
$x=\rho \cos\phi,y=\rho \cos\phi,z=z$. It is assumed that $a>b>0$, and for convenience, that $c>0$ (for a right-handed tube, $c<0$ for a left-handed tube).
The surface lies between the inner and outer cylinders of radius $a-b\le \rho \le a+b$ around the symmetry axis taken here to be the vertical $z$-axis.
Setting $c=0$ reduces this surface to a torus.
By restricting the two angles $\chi$ and $\phi$ to either interval $(-\pi,\pi]$ or $[0,2\pi)$ as convenient, they may be interpreted as coordinate functions on the surface.

The ``central" helical curve corresponds to setting the radius of the vertical circular cross-section of the tube to zero: $b=0$, while $a$ is the radius of the cylinder containing that central helix, and $c$ is its inclination parameter, with ``coiling" angle of inclination $\eta=\arctan(c/a)$. The geometry of this curve parametrized by $\phi$ is easily seen by introducing the unit vectors along the cylindrical coordinate lines
\beq
  \hat\rho = \langle \cos\phi,\sin\phi,0\rangle\,,\
  \hat\phi = \langle -\sin\phi,\cos\phi,0\rangle\,,\
  \hat z = \langle 0,0,1\rangle\,,\
\eeq 
in terms of which it and its tangent have the following representation
\beq
  \vec r\,(\phi) = a\hat\rho + c\phi \hat z \,,\
  \vec r\,'(\phi) = a\hat\phi + c\hat z \,\
\eeq
illustrated in Fig.~\ref{fig:inclination}.

\begin{figure}[p] 
\typeout{*** EPS figure torus1}
\vglue-2cm
\begin{center}
\includegraphics[scale=0.35]{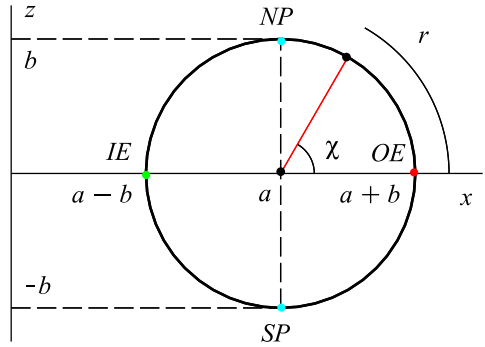}\\[-1cm]
 \includegraphics[scale=0.5]{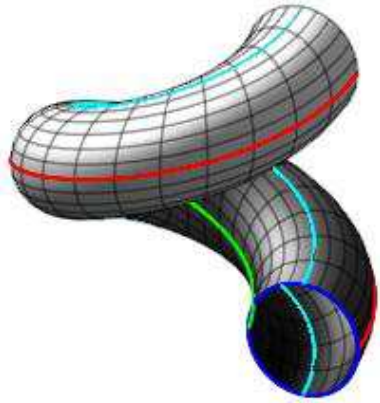}
\end{center}
\vglue-2cm
\caption{
\textbf{Above:}
A vertical half-plane cross-sectional circle of the helical tubular surface built around a helix through the center of this circle whose axis of symmetry is the $z$-axis. This circle in the $x$-$z$ plane is simultaneously rotated around this axis while being translated upwards along that axis ($c>0$), so that the right hand rule wrapping fingers around the helix in the direction in which it is rising (lower figure) puts the thumb up.
\textbf{Below:}
This figure illustrates one turn $0\leq \phi\leq 2\pi$ of the ``unit tube" case $(a,b)=(2,1)$ of a unit circle which is revolved and translated around and along the $z$-axis with inclination angle $\eta=\arctan(2/5)\approx 21.8^\circ$, with an inner equator always a unit distance from the axis, a cavatappo surface model due to Sander Huisman \cite{blog}.
The outer ($\chi=0$, red) and inner ($\chi=\pm\pi$, green) equators are shown together with the ``prime meridian" ($\theta=0$, blue). The northern and southern polar helices ($\chi=\pm \pi/2$, cyan) are equidistant from the two equators.
} 
\label{fig:helix1}
\end{figure}

Fig.~\ref{fig:helix1} illustrates the construction and one complete revolution of a helical tubular surface with its 
inner ($\chi=\pi$) and outer ($\chi=0$) equators marked off.  The grid shown in the computer rendition of the surface consists of the constant $\phi$ circles which result from the intersection of the torus with vertical planes through the symmetry axis (the meridians: $\phi=\phi_0$) and the constant $\chi=\chi_0$ helices (the ``parallels," similar to the parallels in surfaces of revolution). 
The northern ($\chi=\pi/2$) and southern  ($\chi=-\pi/2$) polar helices correspond to the northern and southern polar circles on the torus which in turn generalize the north and south poles on the sphere. The ``radial" arc length coordinate $r=b\,\chi$ and the corresponding angle $\chi$ are measured upwards along the meridians from the outer equator. 
One can also introduce the four hemitubes generalizing the corresponding hemitori of the torus:
the inner ($-\pi/2<\chi< \pi/2$),
the outer ($\pi/2<\chi< 3\pi/2$),
the upper ($0<\chi< \pi$),
and the lower ($-\pi<\chi< 0$)
hemitubes.

If one thinks about how to generate a tubular helical surface with circular cross-sections, there is no reason why not to consider a circular cross-section that is tilted with respect to the vertical direction assumed in the present model (cavatappo 1.0). It is easy to add one additional parameter to this family of basic surfaces which is the angle of tilt of the normal to the plane of the circular cross-section in the $\hat \phi$-$\hat z$ plane of the tangent space up from the horizontal direction. In this enlarged family (cavatappo 2.0), there is an obvious preferred tilt angle which makes the plane of the circular cross-section coincide with the normal plane to the central helical curve. This new orthogonally tilted cavatappo model has very nice mathematical properties which are discussed in a subsequent article \cite{cavatappo20}.

Substituting the differentials of these coordinates into the Euclidean metric $ds^2=dx^2+dy^2+dz^2$ to evaluate the induced metric on the surface, one finds easily for its line element representation
\beq
 ds^2 = \left((a+ b \cos \chi)^2 +c^2\right) d\phi^2 + 2 c\, b\, \cos \chi\, d\phi\, d\chi + b^2 d\chi^2 \,.
\eeq
This metric is independent of $\phi$, which means that it is invariant under translations of the azimuthal angle $\phi$ along the family of parallels, termed helical symmetry \cite{wikiscrew} or more suggestively ``screw-rotation" symmetry. 
The vector field $\partial/\partial \phi$ on the surface which generates these translations is said to be a Killing vector field of this metric. For motion along geodesic curves within the surface, the component of its affinely parametrized geodesic tangent along the Killing vector field remains constant, a ``contant of the motion," or a conserved momentum associated with this symmetry group. This Killing vector field is just the restriction to the surface of the Killing vector field of the Euclidean metric
$\xi_{\rm space} =x\,\partial/\partial y -y\, \partial/\partial x+c\, \partial/\partial z = \partial/\partial \phi + c\, \partial/\partial z$ expressed in either Cartesian or cylindrical coordinates, which generates a rotation about the ``screw axis" by an amount which is proportional to a simultaneous translation along the direction of that axis. In cylindrical coordinates these corkscrew rotations are $\rho\to\rho,\phi \to \phi+t, z\to z+ct$.
Within the surface, these are just translations in $\phi$, so the Killing vector field expressed in the surface coordinates, namely the intrinsic Killing field, is just $\xi=\partial/\partial\phi$.

If we introduce the arclength radial coordinate $r$ and cylindrical radius function $R(r)$ for the helical center curve
\beq
r=b\, \chi\,,\ R=a+b \cos(r/b)\,,
\eeq
the metric can be written in the following form as well as a second form obtained by completing the square on the differential $d\phi$, for which we give two versions in order to highlight the conserved momentum combination we shall see emerging for geodesic motion
\begin{eqnarray}
 ds^2 &=& \left(R^2 +c^2\right) d\phi^2 + 2 c  \cos(r/b)\, d\phi\, dr + dr^2 \nonumber\\
  &=& (R^2 +c^2)\, \left( d\phi+\frac{c\, \cos(r/b)}{R^2+c^2}\, dr \right)^2 + \left(\frac{R^2 +c^2  \sin^2(r/b)}{R^2+c^2}\right)\,dr^2\\ 
  &=&  \frac{ \left((R^2+c^2)d\phi+ c\, \cos(r/b)\,dr \right)^2}{R^2 +c^2}  + \left(\frac{R^2 +c^2  \sin^2(r/b)}{R^2+c^2}\right)\,dr^2
\,. \nonumber
\end{eqnarray}
Completing the square  adapts the metric to the orthogonal decomposition of the tangent space with respect to the intrinsic Killing vector field $\xi=\partial/\partial \phi$. This metric is invariant under reflections through the origin of the coordinates: $(r,\phi) \to (-r,-\phi)$.

Note that the metric determinant (from the second line, an orthogonal frame) is the factor $R^2 +c^2  \sin^2(r/b)$, whose square root provides the surface area density function: $dS=\left(R^2 +c^2  \sin^2(r/b)\right)^{1/2} dr\,d\phi$.
This can be integrated exactly over one revolution of the surface, but the result is an extremely long formula infested with elliptic functions which has little value in being reproduced.

If $(r(\lambda),\phi(\lambda))$ is an affinely parametrized geodesic of this metric on the surface, with tangent
\beq
 U =   \frac{dr}{d\lambda} \partial_r + \frac{d\phi}{d\lambda} \partial_\phi = U^r \partial_r + U^\phi  \partial_\phi
\,,
\eeq
then its orthogonal decomposition is
\beq
 U = U^{\hat r} e_{\hat r} + U^{\hat\phi} e_{\hat\phi}
\eeq
expressed in terms of its components
\begin{eqnarray}
  U^{\hat r} &=& \left(\frac{R^2+c^2 \sin^2(r/b)}{R^2+c^2}\right)^{1/2}  \frac{dr}{d\lambda} \,,\nonumber\\
  U^{\hat \phi} &=& (R^2+c^2)^{1/2}\left(  \frac{d\phi}{d\lambda} + \frac{c \cos(r/b)}{R^2+c^2} \frac{dr}{d\lambda} 
\right ) 
\end{eqnarray}
with respect to the orthonormal frame
\begin{eqnarray}
   e_{\hat r} &=&  \left(\frac{R^2+c^2 \sin^2(r/b)}{R^2+c^2}\right)^{-1/2}
                     \left(\partial_r + \frac{c \cos(r/b)}{R^2+c^2} \partial_\phi\right)\,,
\nonumber\\
   e_{\hat \phi} &=& (R^2+c^2)^{-1/2} \frac{\partial}{\partial \phi}
\end{eqnarray}
whose dual frame is
\beq
  \omega^{\hat r} = \left(\frac{R^2+c^2 \sin^2(r/b)}{R^2+c^2}\right)^{1/2} dr
\,,\
  \omega^{\hat \phi} =  (R^2+c^2)^{1/2} \left(d\phi + \frac{c \cos(r/b)}{R^2+c^2} dr \right) \,.
\eeq

The component of the tangent vector along the Killing vector field is a
conserved screw-angular momentum, i.e., a constant along the geodesic
\beq
  \ell  =  \frac{\partial}{\partial \phi} \cdot U 
    =  \left(R^2 +c^2\right)  \left( \frac{d\phi}{d\lambda}\right) + c\,  \cos(r/b)\, \left( \frac{dr}{d\lambda}\right) 
    =  \left(R^2 +c^2\right)^{1/2} U^{\hat\phi}
\,, 
\eeq
which is exactly the combination occurring in the last form of the metric adapted to the orthogonal decomposition with respect to the radial direction along the meridians.
In an affine parametrization, the square of the length of the tangent is also a constant, which we will call twice the energy
\beq
\left(\frac{ds}{d\lambda}\right)^2
= \frac{\ell ^2}{R^2 +c^2} + \left(\frac{R^2 +c^2  \sin^2(r/b)}{R^2+c^2}\right)\,\left( \frac{dr}{d\lambda}\right)^2
= 2E\,,
\eeq
where we have substituted for $d\phi/d\lambda$ in terms of the screw-angular momentum $\ell$.
The length of the tangent vector in an affine parametrization of a geodesic is a constant along that geodesic, with $E=\frac12$ for an arclength parametrization in which this length is 1.

If we interpret this problem as geodesic motion in the surface in the physics language of motion in space where the affine parameter $\lambda$ plays the role of the time (and $U$ is then called the velocity vector and $U^r,U^\phi$ the velocities), $E$ and $\ell $ are called ``constants of the motion."
Since both $\ell $ and $E$ are constants of the motion, we obtain a single constraint on the square of the ``radial velocity" $dr/d\lambda$,
or equivalently the orthonormal component $U^{\hat r}$, which is of the form
\beq
   {\textstyle\frac12} \left( U^{\hat r}\right)^2 + V = E \,,
\eeq
where
\beq
  V = \frac{\ell ^2}{2(R^2 +c^2) } 
    = \frac{\ell ^2}{2\left((a+b\cos(r/b))^2 +c^2\right) } 
\eeq
acts as an ``effective potential" for the radial motion alone. Note that it is an even function of the radial variable $r$, as well as periodic with period $2\pi b$. Qualitatively, this potential leads to the same kind of radial motion as for the special case $c=0$ of a ring torus thoroughly discussed in the previous article \cite{bobtorus}.
One can fix the potential for $\ell \neq0$ by choosing $\ell $ to have a particular value, and then using the freedom of the affine parametrization to vary the energy $E$ corresponding to energy levels in the graph of the effective potential to describe the allowed interval of radial motion where $V\le E$. When $V=E$ at $r=r_\pm$ corresponding to $dr/d\lambda=0$, turning points of the motion occur. These intervals $r_-\le r\le r_+$ to which such motion is confined are symmetric about $r=0$: $r_-=-r_+$. The outer equator and the corresponding geodesics may be called bound orbits in analogy with the torus problem and more generally 1-dimensional motion in a potential well.
If $E>V$, the motion is unbound with unbounded values of the radial coordinate corresponding to an infinite number of crossings of the inner equator. Local minima of the potential then correspond to stable equilibria, while local maxima correspond to unstable equilibria.
 The outer and inner equators are themselves geodesics which correspond respectively to a stable and unstable equilibrium. 

Note that solving the energy equation for $dr/d\lambda=f(r)$, the problem is reduced to a quadrature
\beq
     \lambda=\int_{r_0}^r f(t)\, dt
\eeq
which must be evaluated numerically,
after which the screw-momentum equation  can be rewritten for $\phi=\phi(r)$ 
\beq
\ell =F(r,dr/d\lambda,d\phi/d\lambda) = F(r,f(r)^{-1},f(r) d\phi/dr)
\eeq
and solved for $d\phi/dr$ and integrated formally  to give $\phi$  as a function of $r$ (only amenable to numerical integration), but this gives no overview of the classification of the orbits which result from the process. The effective potential instead gives a simple visual description of that classification and enables particularly interesting orbits to be studied.

The second order geodesic equations are
\begin{eqnarray}
&&
\mathcal{D} 
{\frac {d^{2}r}{d{\lambda }^{2}}} 
 +  \frac{{c}^{2}}{b}\cos  {\frac {r}{b}} \sin {\frac {r}{b}}  \left( {\frac {dr}{d\lambda }}   \right) ^{2} 
\nonumber\\ &&\qquad
+2\,c\cos  {\frac {r}{b}}   \left( a+b\cos \left( {\frac {r}{b}} \right)  \right) 
   \sin {\frac {r}{b}} {\frac {dr}{d\lambda }}   \frac {d\phi}{d\lambda } 
\nonumber\\ &&\qquad
 +
   \left(\left( a+b\cos \frac {r}{b}  \right) ^{2}+{c}^{2} \right)  
   \left( a+b\cos \frac {r}{b}  \right) \sin \frac {r}{b}   
   \left( \frac {d\phi}{d\lambda} \right) ^{2}
=0\,,
\nonumber\\ &&
\mathcal{D}
\frac {d^{2}\phi}{d{\lambda }^{2}}
  -  \frac{c}{b}\sin \frac {r}{b}   \left( \frac {dr}{d\lambda }   \right) ^{2}
+2 \left( a+b\cos \frac {r}{b}  \right) \sin \frac {r}{b}  \frac {dr}{d\lambda }  \frac {d\phi}{d\lambda }
\nonumber\\ &&\qquad
+c\cos \frac {r}{b}  \left( a+b\cos \frac {r}{b}  \right) \sin \frac {r}{b} \left( \frac {d\phi}{d\lambda}  \right) ^{2}
=0\,,
\nonumber\\ &&
\mathcal{D}=
\left(  \left( a+b\cos \frac {r}{b}  \right) ^{2}+{c}^{2}  \sin^2 \frac {r}{b}  \right) ^{2}
\,.
\end{eqnarray}
Notice that $(dr/d\lambda,d\phi/d\lambda)\to (-dr/d\lambda,-d\phi/d\lambda)$ is a symmetry of this system which maps $\ell \to - \ell $, and so its initial data at a given initial position is reflection symmetric about the origin in the initial tangent space.

Note that $r=0,b\,\pi$ reduces the geodesic equations to $d^2\phi/d\lambda=0$ requiring the azimuthal angle to be linear in the affine parameter $\lambda$ along the inner and outer equators, which are therefore geodesics as in the torus case.

\begin{figure}[p] 
\typeout{*** EPS figure torus2}
\vglue-1cm
\begin{center}
\includegraphics[scale=0.35]{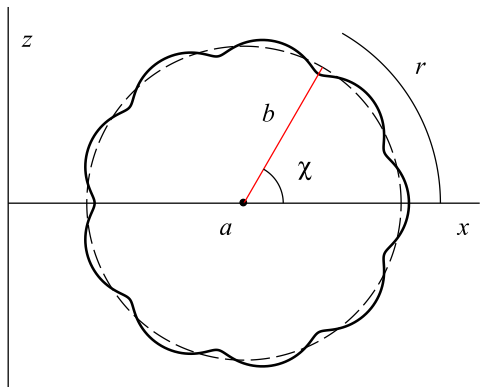}\\[-1cm]
 \includegraphics[scale=0.5]{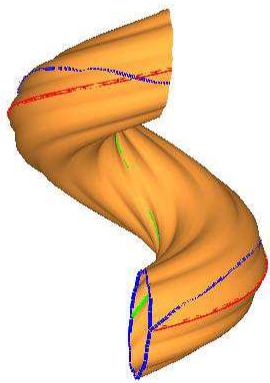}
\end{center}
\caption{
A vertical half-plane cross-sectional circle of the cavatappi (rigati) pasta shape, together with one turn of the surface.
The radial  coordinate $r=b\,\chi$ is no longer an arclength coordinate on the surface, but instead measures the arclength along the original circular cross-section. A bound geodesic is shown which makes a bit more than 1 1/2 oscillations about the outer equator.
} 
\label{fig:helix2}
\end{figure}

If one evaluates the Gaussian curvature function on the surface (half the Riemann scalar curvature) using a computer algebra system, one finds
\beq
K= \frac{1}{b} \frac{\cos(r/b) (a+b\cos(r/b))^3 -b c^2 \sin^4(r/b)}{\left( (a+b\cos(r/b))^2 +c^2 \sin^2(r/b)\right)^2}
\,.
\eeq
The sign of the curvature would be the sign of  $\cos(r/b)$ if it weren't for the explicitly negative term in the numerator, which pushes the zero of curvature from the polar helices towards the outer equator, surrounded by a zone of  positive curvature like the torus.

\section{Generalized helical surfaces: the real cavatappi surface as an instructive example}

One may generalize this family of helical tube surfaces by letting the circular meridian be replaced by any closed plane curve not intersected by the screw-axis, but one needs to be a little more organized about the corresponding calculations. This is essentially an application of the threading splitting approach to a general relativistic spacetime and so serves as a useful example for practice with the associated ideas in a more concrete setting \cite{mfg}.
An instructive example is the actual cavatappi pasta shape, which is ``rigati" (ridged) by a small sinusoidal perturbation of the circular cross-section, with amplitude $d$ and integer frequency $m>1$ (an integer is necessary for periodicity in turn required for the surface to be smooth), leading to $m-1$ outward bulges and $m-1$ inward bulges. This is a 5-parameter family of surfaces described by the equations
\begin{eqnarray}
  x &=& (a+b \cos \chi+d\cos m\chi) \cos \phi\,,\nonumber\\
  y &=& (a+b \cos \chi+d\cos m\chi) \sin \phi\,,\\
  z &=& c\, \phi+b \sin \chi+d\sin m\chi \,,\nonumber
\end{eqnarray}
with $-\infty< \phi , \infty$, $-\phi< \chi \leq \phi$.
The pasta design expert George Legendre \cite{pastadesign} gives the parameter values (only the ratios of the first four are important for the shape, we scale his values down by a factor of 2):
\beq
  a=\frac32\,,\ b=1\,,\ c=\frac45\,, d=\frac{1}{20}\,,\ m=10\,,
\eeq
one complete revolution of which is shown together with its meridian cross-section in Fig.~\ref{fig:helix2}. The pasta shape has 2.5 complete revolutions: $0\le \phi \le 5\pi$. Marked also is a bound geodesic path starting at the outer equator. Actually George uses $c=5/2\pi\approx 0.796 \approx 4/5=0.800$, but the nearby rational value of the parameter slightly simplifies formulas.
 
If one looks closely at the cross-section in Fig.~\ref{fig:helix2}, one sees an asymmetry between the outward bulges and the inward depressions. A symmetrical ridging would have first set $d=0$ and instead replaced $b$ by $b+d\cos(m\chi)$ in Eq.~(\ref{eq:helix}). Clearly  one can study many variations of this theme with the same techniques used below, nor do the particular choice of cavatappo parameter values have any more than an approximate meaning in the real pasta world.

If in addition to $r=b \chi$
we introduce a few abbreviations
\begin{eqnarray}
 R&=&a+b \cos(r/b)+d \cos(mr/b)\,,\nonumber\\
 S&=&b \cos(r/b)+d m\cos(mr/b)\,,\\ 
 T&=&b \sin(r/b)+d m\sin(mr/b)\,,\nonumber
\end{eqnarray}
the metric can be written in the following orthogonal form
\begin{eqnarray}
 ds^2 &=& \left(R^2 +c^2\right) d\phi^2 + 2 \frac{cS}{b}  \, d\phi\, dr 
           + \left(\frac{S^2+T^2}{b^2}\right) dr^2 \nonumber\\
&=& g_{\phi\phi} \, d\phi^2 +2 g_{\phi r} \, d\phi\, dr + g_{rr}\, dr^2\\
&=& M^2 (d\phi+M_r dr)^2 +\gamma_{rr} \, dr^2
\,,\nonumber\\
&=& (\omega^{\hat\phi})^2 +(\omega^{\hat r})^2 \,,
\end{eqnarray}
which requires computation of the following quantities:
\begin{eqnarray}
  M &=&(g_{\phi\phi})^{1/2} =\left(R^2 +c^2\right)^{1/2}
\,,\ 
    M_r = \frac{g_{\phi r}}{g_{\phi\phi}} =  \frac{cS}{b (R^2 +c^2)}  
\,,\nonumber\\
  \gamma_{rr} &=& g_{rr} -M^2 M_r^2  = \frac{S^2+T^2-c^2 S^2 /(R^2+c^2)}{b^2}
\,.
\end{eqnarray}
This corresponds to the orthogonal frame and dual frame
\begin{eqnarray}
&& e_r = \partial_r-M_r\, \partial_\phi\,,\
   e_\phi = \partial_\phi\nonumber\\
&& \omega^r = dr\,,\
   \omega^\phi = d\phi+M_r\, dr
\end{eqnarray}
and orthonormal frame and dual frame
\begin{eqnarray}
&& e_{\hat r} = \gamma_{rr}^{-1/2} e_r\,,\
   e_{\hat\phi} = M^{-1/2} e_\phi\nonumber\\
&& \omega^r = \gamma_{rr}^{1/2} \omega^r\,,\
   \omega^\phi = M \omega^\phi\,.
\end{eqnarray}
Note that $M \gamma_{rr}^{1/2} = (R^2+c^2)(S^2+T^2)- c^2 S^2$ 
is the determinant of the metric and determines the surface volume element $M \gamma_{rr}^{1/2} dr\, d\phi$ of the tubular surface.

\begin{figure}[p] 
\typeout{*** EPS figures IVP}
\vglue-1cm
\begin{center}
\includegraphics[scale=0.3]{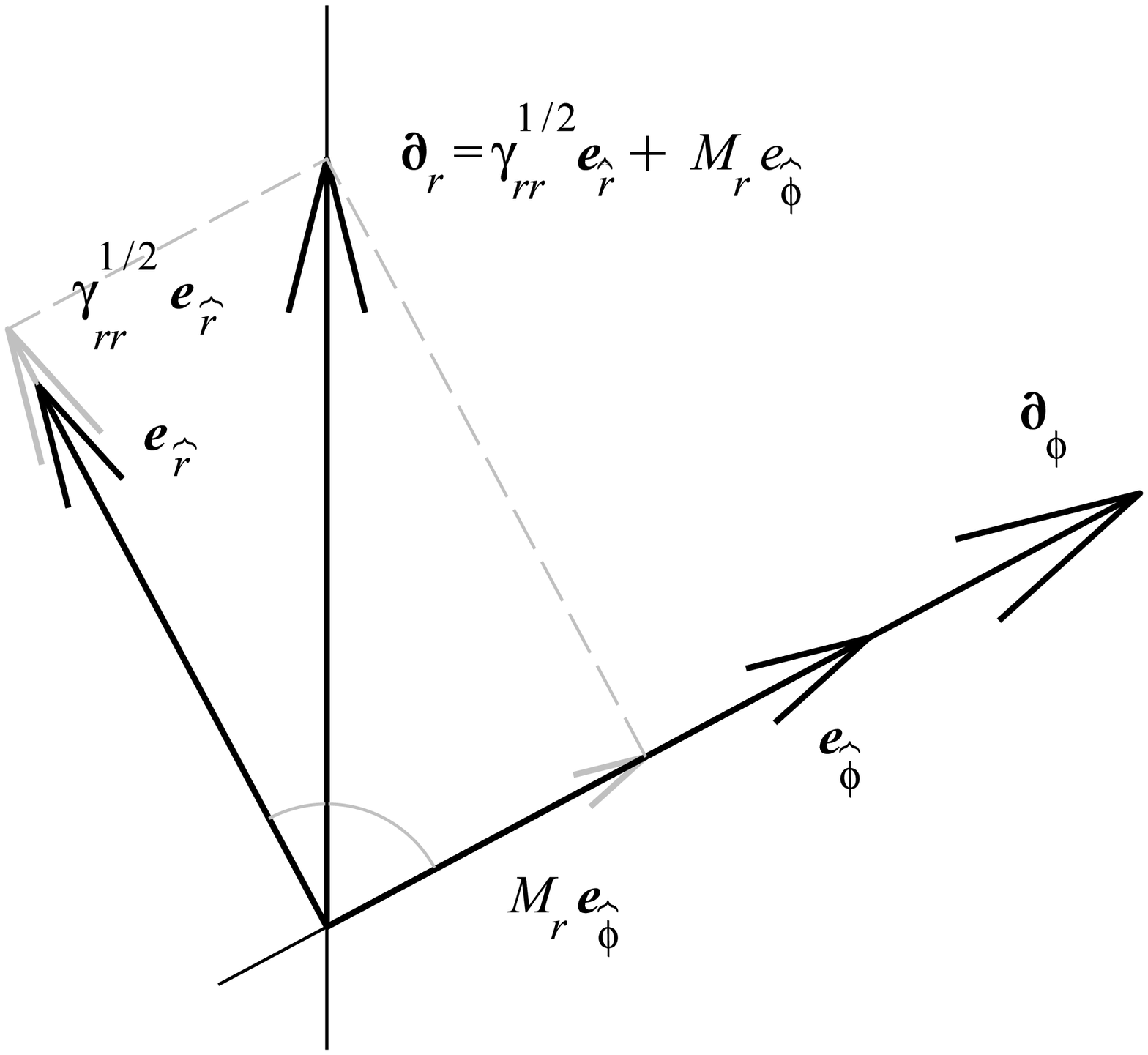}\\
 \includegraphics[scale=0.3]{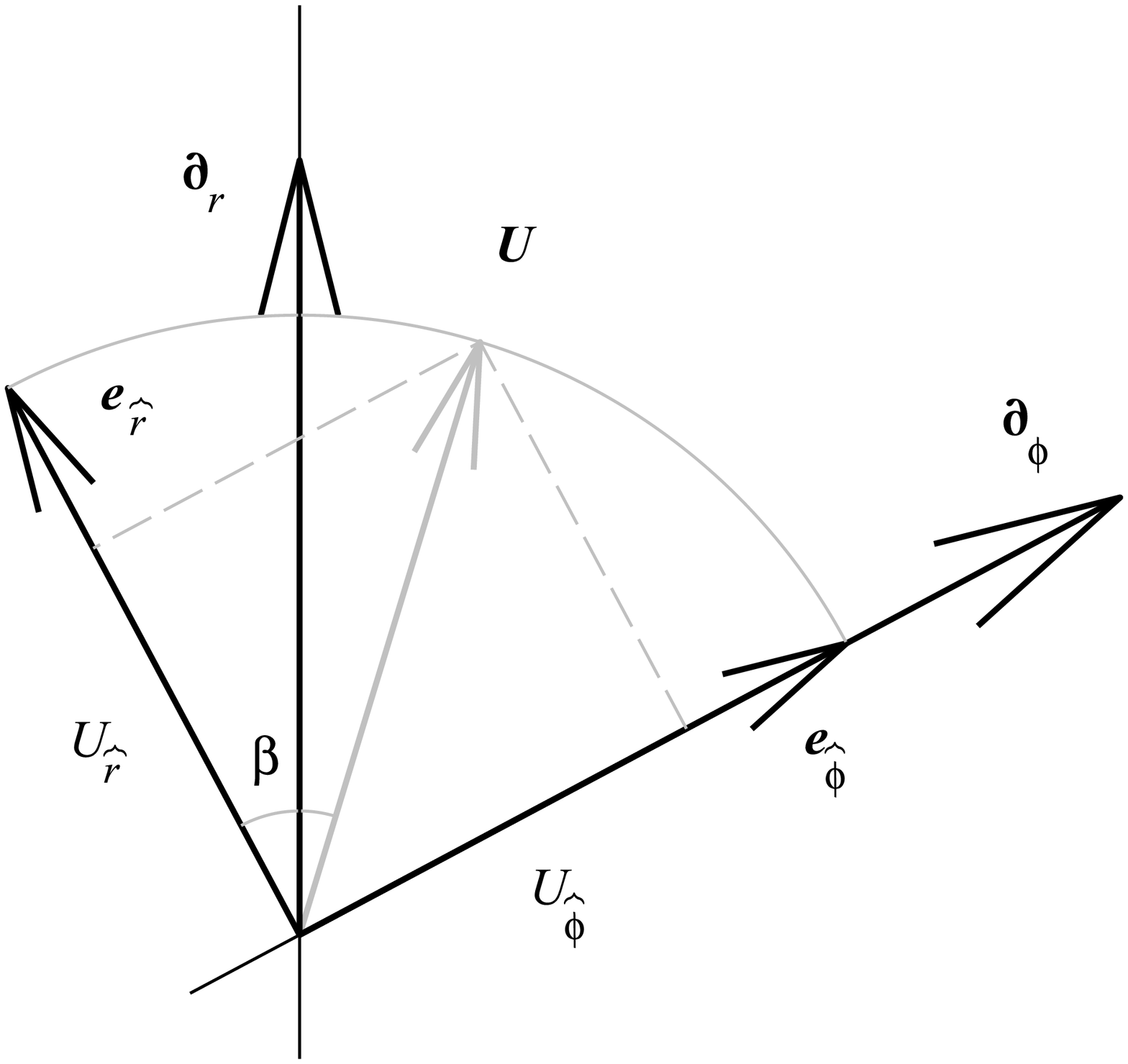}
\end{center}
\caption{
The orthogonal decomposition of the radial direction with respect to the parallels in the tangent space (top). The initial data is specified as a unit vector with angle $\beta$ measured down from the direction orthogonal to the parallels (bottom). 
} 
\label{fig:ivp}
\end{figure}

If $(r(\lambda),\phi(\lambda))$ is an affinely parametrized geodesic of this metric on the surface, its tangent is
\beq
 U =   \frac{dr}{d\lambda} \partial_r + \frac{d\phi}{d\lambda} \partial_\phi
  = U^{\hat r} e_{\hat r} + U^{\hat\phi} e_{\hat\phi}
\eeq
and has orthonormal  components
\begin{eqnarray}
  U^{\hat r} = \gamma_{rr}^{1/2}  \frac{dr}{d\lambda} \,,
\quad 
  U^{\hat \phi} = M \left(  \frac{d\phi}{d\lambda} + M_r \frac{dr}{d\lambda}  
\right ) \,.
\end{eqnarray}
These can be inverted to yield
\begin{eqnarray}  \label{eq:velocityrelation}
  U^r= \frac{dr}{d\lambda} = \gamma_{rr}^{-1/2}  U^{\hat r} \,,
\quad 
 U^\phi = \frac{d\phi}{d\lambda} = M^{-1} U^{\hat \phi} - M_r \gamma_{rr}^{-1/2}   U^{\hat r}\,. 
\end{eqnarray}

\begin{figure}[t] 
\typeout{*** EPS figure helix potential 2}
\vglue-1cm
\begin{center}
\includegraphics[scale=0.35]{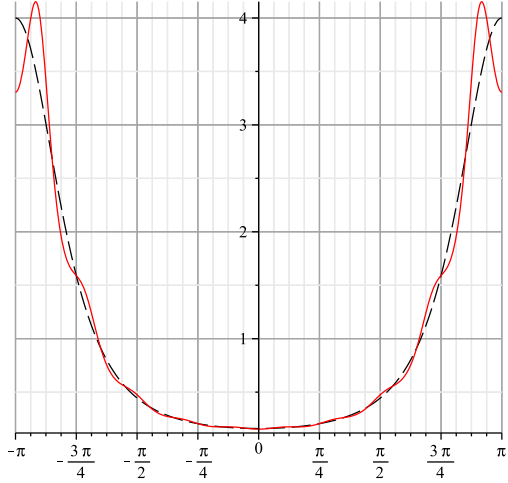}
\end{center}
\caption{
The effective potential for geodesic motion on the cavatappi pasta surface (solid) and its smoothed version (dashed). Since the frequency $m=10$ is even, there are an odd number of bulges, so in addition to the the inner equator geodesic, there is a pair of unstable geodesics on either side of it on the two nearby innermost bulges closest to the screw axis.
} 
\label{fig:helixpotential2}
\end{figure}

The component of the tangent vector along the Killing vector field  (the conserved screw-angular momentum) is a constant along the geodesic
\beq
  \ell  =  M U^{\hat \phi} = M^2 \left( \frac{d\phi}{d\lambda}+M_r \frac{dr}{d\lambda} \right)
\,, 
\eeq
as is its length, half of which we call the energy
\beq
   {\frac12}\left( U^{\hat \phi} \right)^2 +   {\frac12}\left( U^{\hat r} \right)^2
 =  {\frac12}\left( U^{\hat r} \right)^2 +  {\frac12}\frac{\ell ^2}{M^2} =  E
\eeq
so again
\beq
   V =  {\frac12}\frac{\ell ^2}{M^2}
     =  {\frac12}\frac{\ell ^2}{(a+b \cos(r/b)+d \cos(mr/b))^2+c^2}
\eeq
acts as an effective potential for the radial motion. If the screw-angular momentum is nonzero, then we might as well set it equal to 1 and use the energy parameter to distinguish initial data.
Fig.~\ref{fig:helixpotential2} shows the nonzero screw-angular momentum potential for the cavatappi pasta surface together with the corresponding smooth shape with $d=0$.  

When $\ell \neq0$, of course we can solve the energy constraint on the radial velocity for that derivative $dr/d\lambda=f(r)$ and then integrate $d\lambda = dr/f(r)$, but this is extremely complicated and must be done numerically.
If $\ell =0$, then we get a simple first order differential equation to integrate for the angle $\phi$ as a function of the radial variable $r$ 
\begin{eqnarray}
   && \frac{d\phi}{d\lambda}+M_r \frac{dr}{d\lambda}=0 \to 
  d\phi =-M_r dr \to \nonumber\\
   && \phi =-\int_0^r M_r (w)\, dw = -\int_0^r \frac{c}{b}\left(b \cos(w/b)+d m\cos((mw/b)\right)dw
\,.
\end{eqnarray}
In particular one loop of the surface corresponds to 
\beq
 \Delta \phi =-\int_0^{2\pi} M_r (w)\, dw \neq0\,.
\eeq
which does not vanish since this is a function which is even about the inner equator $\phi=\pi$ which is more negative near the inner equator than it is near the outer equator $\phi=0$. Thus these are not closed geodesics.

If we introduce polar coordinates in the tangent plane (``velocity space") adapted to the orthonormal frame, then the velocity has constant length $(2E)^{1/2}$ because of the affine parametrization assumed for the geodesic, so we get the representation
\beq
        \langle U^{\hat r}, U^{\hat\phi}\rangle = (2E)^{1/2} \langle \cos\beta,\sin\beta \rangle
\eeq
in terms of the angle measured from the direction orthogonal to the parallels, which in turn implies the constancy of 
\beq\label{eq:Msinbeta}
  \frac{\ell}{(2E)^{1/2}} = M\, \sin\beta\,.
\eeq
For a given energy that sets the affine scale, and for a given initial position on the surface, the initial angle $\beta$ parametrizes the family of geodesics which emanate from this point. If we choose to numerically solve the coordinate component form of the geodesic equations rather than the orthonormal component form, then  Eq.~(\ref{eq:velocityrelation}) expressed in terms of $\beta$ gives us geometric initial velocity data for those second order differential equations. Setting $E=1/2$ then corresponds to an arclength parametrization. If instead we set $\ell=1$, then the relation (\ref{eq:Msinbeta}) sets the energy level according to the initial values of $r$ and $\beta$.

At this point one can study the details of this system in the same approach taken in the preceding article \cite{bobtorus} written to be understood by undergraduates who have already taken multivariable calculus, differential equations and linear algebra. A computer algebra system is essential to make concrete calculations. While George Legendre chose Mathematica, the author is a long time Maple user.

\section{Italy and Tensor Calculus/Differential Geometry}

Dirk Stuik also recalled that when asked what he liked best about Italy, Einstein responded ``spaghetti and Levi-Civita."\\
---quote from a talk at a 1996 AMS Meeting \cite{spaghetti}.
\vskip\baselineskip

While there is widespread awareness of the association between pasta and Italy, few people realize the key role Italian mathematicians played in developing the tools of tensor analysis and differential geometry which made Einstein's theory of general relativity possible. While the initial ideas of curvature of surfaces were developed by Gauss, and Riemann generalized them to arbitrary dimension, it was Gregorio  Ricci-Curbastro who provided the mathematical infrastructure of tensor calculus to make it usable for calculating, and his student/collaborator Tullio Levi-Civita \cite{rossana} joined this enterprise, codifying the ideas into a landmark article in 1900 known for a time as Ricci Calculus \cite{ricci}. 
Einstein was not very mathematical, and it was this article that his friend Marcel Grossmann \cite{grossmann} used in 1912--1914 to teach him enough to get serious about Einstein's attempts to generalize special relativity, and Levi-Civita participated in Einstein's progress through and after the birth of the final theory of general relativity published in 1916, studying its solutions and introducing the key idea of parallel transport in 1917 \cite{paralleltransport}. This latter concept plays a fundamental role in both using general relativity and in understanding differential geometry. Levi-Civita presented this theory in lecture notes put together with Enrico Persico (friend of Enrico Fermi, and editor of his collected works) in 1925 and translated into English in 1927 as \textit{The Absolute Differential Calculus} \cite{levi,wiki}.

The present family of surfaces provide an interesting playground for students to learn a bit of geodesic theory while perhaps having some fun. Yet unexplored is the extrinsic curvature properties of these surfaces. However, the discussion so far also shows how essential computer algebra systems are for actually doing the calculations once some intelligent framework is imposed on the problem.

\section*{Acknowledgements}

This work would not have been possible without the initial Maple torus worksheet shared with me by my colleague Klaus Volpert at just the right moment.

\end{document}